\documentclass{amsproc}

\usepackage{graphicx}

\usepackage{amscd}

\usepackage{amsmath}

\usepackage{amsfonts}

\usepackage{amssymb}

\theoremstyle{plain}

\newtheorem{theorem}{Theorem}[section]

\newtheorem{corollary}[theorem]{Corollary}

\newtheorem{definition}[theorem]{Definition}

\newtheorem{example}[theorem]{Example}

\newtheorem{examples}[theorem]{Examples}

\newtheorem{proposition}[theorem]{Proposition}

\numberwithin{equation}{section}

\begin{document}

\title[Branched Coverings]{$C^{\ast}$-algebras Associated with Branched Coverings}

\author{Valentin Deaconu}

\address{Department of Mathematics, University of Nevada, Reno, NV 89557}

\email{vdeaconu@math.unr.edu}

\urladdr{http://unr.edu/homepage/vdeaconu}

\author{Paul S. Muhly}

\address{Department of Mathematics, University of Iowa, Iowa City, IA 52242}

\email{muhly@math.uiowa.edu}

\thanks{Research supported in part by grants from the National Science Foundation.}

\date{June 25, 1999}



\begin{abstract}In this note we analyze the $C^{\ast}$-algebra associated with a branched
covering both as a groupoid $C^{\ast}$-algebra and as a Cuntz-Pimsner algebra.
We determine conditions when the algebra is simple and purely infinite. We
indicate how to compute the K-theory of several examples, including one related to rational maps on the Riemann sphere.

\end{abstract}

\maketitle

\section{Introduction}

Given a branched covering $\sigma:X\rightarrow X$ of a locally compact space
$X$, we define its $C^{\ast}$-algebra to be the $C^{\ast}$-algebra of the
r-discrete groupoid $\Gamma$ associated by Renault to the corresponding
partially defined local homeomorphism $T$. More precisely,
\[
\Gamma=\Gamma(X,\sigma)=\{(x,m-n,y)\mid m,n\in\mathbb{N},x\in dom(T^{m}),y\in
dom(T^{n}),T^{m}x=T^{n}y\},
\]
where $T$ is the restriction of $\sigma$ to the nonsingular set $X\setminus
S=U$. Here $dom (T^k)$ is the domain of $T^k$.

It turns out that this $C^{\ast}$-algebra is isomorphic to an augmented
Cuntz-Pimsner algebra associated with a $C^{\ast}$-correspondence $(A,E)$,
where $A=C_{0}(X),E=\overline{C_{c}(U)}$, viewed as an $A$-Hilbert module via
the right multiplication
\[
(\xi f)(x)=\xi(x)f(\sigma(x)),\xi\in E,f\in A,x\in U,
\]
and inner product
\[
\langle\xi,\eta\rangle(x)=\sum_{\sigma(y)=x}\overline{\xi(y)}\eta(y).
\]
The left multiplication is given by the map
\[
\varphi:A\rightarrow L(E),\;(\varphi(f)\xi)(x)=f(x)\xi(x).
\]
If the singular set is empty, so that $U=X$, we recover previous results (see
\cite{De4}, \cite{Del}).

We consider several examples of $C^{\ast}$-algebras arising from branched
coverings, and indicate how to compute their K-theory. These examples include some of the
algebras considered by R.Exel in \cite{Ex}, in the case of a partial homeomorphism.

{\bf Acknowledgements}. The first author would like to thank Alex Kumjian and Bruce Blackadar for helpful discussions.

\section{Branched coverings and groupoids}

We collect here some facts about branched coverings and some examples for
future reference.

\begin{definition}
Let $X,X^{\prime}$ be locally compact, second countable Hausdorff spaces and
let $S\subset X,\;S^{\prime}\subset X^{\prime}$ be closed subsets such that
$U=X\setminus S$ and $V=X^{\prime}\setminus S^{\prime}$ are dense in $X$ and
$X^{\prime}$, respectively. A continuous surjective map $\sigma:X\rightarrow
X^{\prime}$ is said to be a \emph{branched covering} with \emph{branch sets}
$S$ (upstairs) and $S^{\prime}$ (downstairs) if

\begin{enumerate}
\item  the components of preimages of open sets of $X^{\prime}$ are a basis
for the topology of $X$ (in particular $\sigma$ is an open map),
\item $\sigma(S)=S^{\prime},\,\,\sigma(U)=V$, and
\item $\sigma\mid_{U}$ is a local homeomorphism.
\end{enumerate}
\end{definition}

We mention that in the original definition given by Fox (see \cite{Fo}),
$X\setminus S$ and $X^{\prime}\setminus S^{\prime}$ are supposed to be
connected and $S$ and $S^{\prime}$ are supposed to be of codimension 2 for
topological reasons. We make no such restrictions here; in particular, we
allow disconnected $X\setminus S$ and $X^{\prime}\setminus S^{\prime}$.

\begin{examples}
\begin{enumerate}
\item [ a)](Folding the interval) Let $X=X^{\prime}=[0,1]$, and the map
\[
\sigma(t)=\left\{
\begin{array}
[c]{ll}%
2t,\;0\leq t\leq1/2 & \\
2-2t,\;1/2\leq t\leq1. &
\end{array}
\right.
\]
Then $S=\{1/2\},\,\,S^{\prime}=\{1\}.$

\item[ b)] Let $X=X^{\prime}=\mathbf{D}^{2}$ (the closed unit disc),
$\sigma(z)=z^{p}\,\,(p\geq2).$ Then $S=S^{\prime}=\{0\}.$

\item[ c)] Let $X=\mathbf{S}^{k}$ be the $k$-sphere, \thinspace$X^{\prime
}=\mathbf{D}^{k}$ be the $k$-disc,  and let $\sigma$ be the
projection onto the equatorial plane. Then $S=S^{\prime}=\mathbf{S}^{k-1}$
(identified with the equator).

\item[ d)] Let $G$ be a finite group acting non-freely on a compact manifold
$X$, let $X^{\prime}=X/G$ be the corresponding orbifold, and let
$\sigma:X\rightarrow X^{\prime}$ be the quotient map. Then $\sigma$ is a
branched covering with $S=$ the set of points $s$ with nontrivial isotropy
group $G_{s}$. The previous examples are just particular cases of this, but
note that not every branched covering is associated with a group action.

\item[ e)] Let $X$ be the unit circle, identified with the one-point
compactification of $\mathbb{R}$. Consider $\sigma:\mathbb{R}\rightarrow X$
the ususal map which wraps $\mathbb{R}$ around the circle. Here $\sigma$ is not defined everywhere, but the groupoid construction will make sense.

\item[ f)] Let $p:\mathbb{C}\rightarrow\mathbb{C}$ be a polynomial map of
degree $k\geq2$ , and define $p(\infty)=\infty$. We get a map $p:\mathbf{S}%
^{2}\rightarrow\mathbf{S}^{2}$, which is a branched covering with
\[
S^{\prime}=\{w\in\mathbf{S}^{2}\mid\mbox{ the equation }\;p(z)=w\mbox
{ has multiple roots}\},
\]
$S=p^{-1}(S^{\prime})$. In a similar way, each rational map $q:\mathbf{S}%
^{2}\rightarrow\mathbf{S}^{2}$ gives a branched covering.

\item[ g)] If $\sigma:U\rightarrow V$ is a homeomorphism between two open
subsets of $X$, then we are in the situation considered by Exel in \cite{Ex},
by taking $A=C_{0}(X),I=C_{0}(V),J=C_{0}(U)$.

\end{enumerate}
\end{examples}

While one can build $C^{\ast}$-correspondences from any branched covering
$\sigma:X\rightarrow X^{\prime}$, in this note we shall restrict our attention
to spaces covering themselves in a branched fashion.

For a locally compact space $X$ and a local homeomorphism $T$ from an open
subset $dom(T)$ of $X$ onto an open subset $ran(T)$ of $X$, Renault denotes by
$Germ(X,T)$ the groupoid of germs of $\mathcal{G}(X,T)$, which in turn is the
full pseudogroup generated by restrictions $T\mid_{Y}$, where $Y$ is an open
subset of $X$ on which $T$ is injective. He proves that the groupoid of germs
coincides with the semidirect product groupoid
\[
\Gamma(X,T)=\{(x,m-n,y)\;\mid\;m,n\in\mathbb{N},x\in dom(T^{m}),y\in
dom(T^{n}),T^{m}x=T^{n}y\},
\]
if $(X,T)$ is essentially free in the sense that for all $m,n$, there is no
nonempty open set on which $T^{m}$ and $T^{n}$ agree.

The topology on $\Gamma$ is generated by the open sets
\[
{}(Y;m,n;Z)=\{(y,m-n,z)\;\mid\;(y,z)\in Y\times Z,T^{m}(y)=T^{n}(z)\},
\]
where $Y$ and $Z$ are open subsets on which $T^{m}$ and $T^{n}$ are injective,
respectively. Since the range and the source maps are local homeomorphisms,
$\Gamma$ becomes a locally compact (Hausdorff) r-discrete groupoid, and so we
may consider its $C^{\ast}$-algebra.

Denote by $c:\Gamma(X,T)\rightarrow\mathbb{Z}$ the cocycle defined by
$c(x,k,y)=k$, i.e., $c$ is the so-called \emph{position cocycle}. Then
$c^{-1}(0)$ is an equivalence relation $R(X,T)$, which is the increasing union
of sub-relations
\[
R_{N}=\{(x,y)\in X\times X\;\mid\;\exists\;n\leq N\;\;{\mbox{with}}\;x,y\in
dom(T^{n}),T^{n}x=T^{n}y\}.
\]

Notice that $R_{0}$ is the diagonal of $X\times X$, and each $R_{N}$ is
r-discrete, because the range and source maps are local homeomorphisms.

\section{$C^{\ast}$-algebras defined by a branched covering}

\begin{definition}
Given a branched covering $\sigma:X\rightarrow X$ with branch set $S$,
consider the local homeomorphism $T:X\setminus S\rightarrow X\setminus
\sigma(S)$, and assume that $(X,T)$ is essentially free. We define $C^{\ast
}(X,\sigma)$ to be $C^{\ast}(\Gamma(X,T))$.
\end{definition}

Let $\displaystyle
A=C_{0}(X)$, and let $E=\overline{C_{c}(U)}$, where $U=X\setminus S$, with the
structure of a Hilbert $A$-module given by the formulae:
\[
(\xi f)(x)=\xi(x)f(\sigma(x)),\;\xi\in E,f\in A,x\in U,
\]%
\[
\langle\xi,\eta\rangle(y):=\sum_{\sigma(x)=y}\overline{\xi(x)}\eta
(x),y\in\sigma(U),\;\xi,\eta\in E.
\]

In other words, the inner product is given by $\langle\xi,\eta\rangle
=P(\bar{\xi}\eta)$, where $P$ is the extension of $P:C_{c}(U)\rightarrow
C_{c}(\sigma(U)),$
\[
(P\xi)(y)=\sum_{\sigma(x)=y}\xi(x).
\]
Note that the inner products generate the ideal $C_{0}(\sigma(U))$ in $A$. The
left module structure on $E$ is defined by the equation
\[
\varphi:A\rightarrow L(E),\;(\varphi(f)\xi)(x)=f(x)\xi(x)\;f\in A,\xi\in E.
\]
It is straightforward to verify that $\varphi(f)$ is in $L(E)$ with adjoint
$\varphi(\bar{f}),f\in A.$ It is also straightforward to verify that $\varphi$
is injective. Thus, we may form the augmented Cuntz-Pimsner algebra
$\tilde\mathcal{O}_{E}$.

\begin{theorem}
\label{isomorphism}The $C^{\ast}$-algebras $\tilde\mathcal{O}_{E}$ and $C^{\ast
}(\Gamma(X,T))$ are isomorphic.
\end{theorem}

The proof of the theorem will be given in several steps, using the notion of
Fell bundle.

We recall briefly the Pimsner construction from \cite{Pi}. A $C^{\ast}%
$-\emph{correspondence }is a pair $(E,A)$, where $E$ is a (right) Hilbert
module over a $C^{\ast}$-algebra $A$, and where $A$ acts to the left on $E$
via a $\ast$-homomorphism $\varphi:A\rightarrow L(E)$, from $A$ to the bounded
adjointable module maps on $E$. We shall always assume that our map $\varphi$
is injective. The module $E$ is not necessarily full, in the sense that the
span of the inner products $\langle E,E\rangle$ may be a proper ideal of $A$.
Given a $C^{\ast}$-correspondence $(E,A)$, Pimsner constructs a $C^{\ast}
$-algebra $\mathcal{O}_{E}$, which generalizes both the crossed products by
$\mathbb{Z}$ and the Cuntz-Krieger algebras. The C*-algebra generated by $\mathcal{O}_E$ and $A$ is denoted by $\tilde\mathcal{O}_E$, and is called  the
{\em augmented} \emph{Cuntz-Pimsner algebra} of the correspondence. The algebra
$\mathcal{O} _{E}$ is a quotient of the generalized Toeplitz algebra
$\mathcal{T}_{E}$ generated by the creation operators $\mathcal{T}_{\xi
},\,\,\xi\in E$ on the Fock space $\displaystyle\mathcal{E}_{+}=\bigoplus
_{n=0}^{\infty}E^{\otimes n}$. Here $E^{\otimes0}=A$, and for $n\geq1,$
$E^{\otimes n}$ denotes the $n$-th tensor power of $E$, balanced via the map
$\varphi$. The creation operators $\mathcal{T}_{\xi}$, $\xi\in E$, are defined
by the formulae $\mathcal{T}_{\xi}a=\xi a,$ for $a\in A$, and $\mathcal{T}%
_{\xi}(\xi_{1}\otimes...\otimes\xi_{n})=\xi\otimes\xi_{1}\otimes...\otimes
\xi_{n},$ for $\xi_{1}\otimes...\otimes\xi_{n}\in E^{\otimes n}$.

A $C^{\ast}$-correspondence should be viewed as a generalization of an
endomorphism of a $C^{\ast}$-algebra. Just as an endomorphism of a $C^{\ast}
$-algebra can be ``extended'' to an automorphism of a larger algebra, so a
general $C^{\ast}$-correspondence can be ``extended'' to an ``invertible''
correspondence. We will not go into much detail here except to indicate how
this leads to another presentation of $\mathcal{O}_{E}$.

Pimsner considers a new pair $(E_{\infty},\mathcal{F}_{E})$, where
$\mathcal{F}_{E}$ is the $C^{\ast}$-algebra generated by all the compact
operators $K(E^{\otimes n}),\;n\geq0$ in $\displaystyle\lim_{\longrightarrow
}L(E^{\otimes n})$, and $E_{\infty}=E\otimes\mathcal{F}_{E}$. The advantage is
that $E_{\infty}$ becomes an $\mathcal{F}_{E}$-$\mathcal{F}_{E}$ bimodule,
such that the dual or adjoint module $E_{\infty}^{\ast}$ is also an
$\mathcal{F}_{E}$-$\mathcal{F}_{E}$ bimodule. More acurately, $E_{\infty}$ has
two $\mathcal{F}_{E}$-valued inner products with respect to which $E_{\infty}$
satisfies all the axioms of an imprimitivity bimodule over $\mathcal{F}_{E}$,
except possibly the one asserting that the left inner product is full. The
$C^{\ast}$-algebra $\mathcal{O}_{E}$ is represented on the two-sided Fock
space
\[
\mathcal{E}_{\infty}=\bigoplus_{n\in\mathbb{Z}}E_{\infty}^{\otimes n},
\]
where for $n<0,E_{\infty}^{\otimes n}$ means $(E_{\infty}^{\ast})^{\otimes-n}%
$. It is isomorphic to the $C^{\ast}$-algebra generated by the multiplication
operators $M_{\xi}\in L(\mathcal{E}_{\infty})$, where for $\xi\in E_{\infty}$,
$M_{\xi}\omega=\xi\otimes\omega$.

Given a branched covering $\sigma:X\rightarrow X$, we first want to identify
the $C^{\ast}$-algebra $\mathcal{F}_{E}$. Recall that $T:U\rightarrow
\sigma(U)$ is the local homeomorphism associated to $\sigma$.

Note that $E\otimes_{\varphi}E$ is a quotient of $\overline{C_{c}(U)}%
\otimes\overline{C_{c}(U)}$, where we identify $\xi f\otimes\eta$ with
$\xi\otimes\varphi(f)\eta$ for any $\xi,\eta\in E$ and any $f\in A$. Therefore
$E\otimes_{\varphi}E$ can be identified, as a vector space, with the
completion of the compactly supported continuous functions on the set
$U\cap\sigma(U)$. In a similar way, $E^{\otimes n}$ is identified (as a vector
space) with $\overline{C_{c}(U\cap\sigma(U)\cap...\cap\sigma^{n-1}(U))}$. The
structure of a Hilbert $A$-module on $E^{\otimes n}$ is given by the
equations
\[
(\xi f)(x)=\xi(x)f(\sigma^{n}(x)),
\]
and
\[
\langle\xi,\eta\rangle_{n}=P_{n}(\bar{\xi}\eta).
\]
Here $P_{n}:C_{c}(U\cap\sigma(U)\cap... \cap\sigma^{n-1}(U))\rightarrow
C_{c}(\sigma^{n}(U))$ is given by
\[
P_{n}(\xi)(y)=\sum_{\sigma^{n}(x)=y}\xi(x).
\]

\begin{proposition}
\label{Prop 3.1}The $C^{\ast}$-algebra $K(E)$ is isomorphic to $C^{\ast
}(R(T))$, where
\[
R(T)=\{(x,y)\in U\times U\mid\;T(x)=T(y)\}
\]
is the equivalence relation associated with $T=\sigma\mid_{U}$.

More generally, $K(E^{\otimes n})\simeq C^{\ast}(R(T^{n}))$, where
\[
R(T^{n})=\{(x,y)\in dom(T^{n})\times dom(T^{n})\mid\;T^{n}(x)=T^{n}(y)\}
\]
is the equivalence relation associated with $T^{n}$.
\end{proposition}

\begin{proof}
It is known that $K(E)=E\otimes E^{\ast}$, the tensor product balanced over
$A$, where $E^{\ast}$ is the adjoint of $E$. Since $\xi f\otimes\eta^{\ast
}=\xi\otimes f\eta^{\ast}$, it follows that, as a set, $K(E)=\overline
{C_{c}(R(T))}$. The multiplication of compact operators is exactly the
convolution product on $C_{c}(R(T))$, therefore, as $C^{\ast}$-algebras,
$K(E)=C^{\ast}(R(T))$.

In the same way, using the fact that $K(E^{\otimes n})=(E^{\otimes n}%
)\otimes(E^{\otimes n})^{\ast}$, we get $K(E^{\otimes n})=C^{\ast}(R(T^{n}))$.
\end{proof}

If we take
\[
R_{N}:=\bigcup_{n=0}^{N} R(T^{n}),
\]
where $R_{0}$ is the diagonal of $X$, then the natural inclusion $C^{*}%
(R_{N})\rightarrow C^{*}(R_{N+1})$ is induced by the map
\[
L(E^{\otimes N})\rightarrow L(E^{\otimes N+1}), \;\; F\mapsto F\otimes I.
\]

\begin{corollary}
We have $\displaystyle\mathcal{F}_{E}=\lim_{\longrightarrow}C^{\ast}(R_{N}).$
In particular, $\mathcal{F}_{E}$ is isomorpic to $C^{\ast}(R(X,T))$, where
$R(X,T)$ is defined at the end of the previous section.
\end{corollary}

In order to establish an isomorphism between $C^{\ast}(\Gamma)$ and
$\tilde\mathcal{O}_{E}$ in our setting, we follow \cite{De5} and use the notion of
Fell bundle. We show that both $C^{\ast}(\Gamma)$ and $\tilde\mathcal{O}_{E}$ are
isomorphic to the $C^{\ast}$-algebra associated to isomorphic Fell bundles
over $\mathbb{Z}$. This point of view was suggested by Abadie, Eilers and Exel
in \cite{AEE}. We recall the definition of a Fell bundle and of the associated
$C^{\ast}$-algebra (cf. \cite{Ku2} for a more general situation).

\begin{definition}
Consider a Banach bundle $p:$\emph{$\mathcal{B}$}$\rightarrow\mathbb{Z}$. A
\emph{multiplication} on \emph{$\mathcal{B}$ }is a continuous map
\emph{$\mathcal{B}\times\mathcal{B}\rightarrow\mathcal{B}$ }($(b_{1}%
,b_{2})\rightarrow b_{1}b_{2}$) which satisfies the conditions

\begin{enumerate}
\item [ a)]$p(b_{1}b_{2})=p(b_{1})+p(b_{2}),b_{1},b_{2}\in$\emph{$\mathcal{B}%
$, }
\item[ b)] \emph{$\mathcal{B}_{k}\times\mathcal{B}_{l}\rightarrow
\mathcal{B}_{k+l}$ }is bilinear,

\item[ c)] $(b_{1}b_{2})b_{3}=b_{1}(b_{2}b_{3}),$

\item[ d)] $\left\|  b_{1}b_{2}\right\|  \leq\left\|  b_{1}\right\|
\,\left\|  b_{2}\right\|  $.
\end{enumerate}

An \emph{involution} is a continuous map \emph{$\mathcal{B}\rightarrow
\mathcal{B},\;$}$b\mapsto b^{\ast}$\emph{, }which satisfies

\begin{enumerate}
\item [ e)]$p(b^{\ast})=-p(b),b\in$\emph{$\mathcal{B},$}

\item[ f)] \emph{$\mathcal{B}_{k}\rightarrow\mathcal{B}_{-k}$ }is conjugate
linear $\forall k\in\mathbb{Z}$,

\item[ g)] $b^{\ast\ast}=b$.
\end{enumerate}

The bundle \emph{$\mathcal{B}$ }together with these maps is said to be a
\emph{Fell bundle} if in addition

\begin{enumerate}
\item [ h)]$(b_{1}b_{2})^{\ast}=b_{2}^{\ast}b_{1}^{\ast}$,

\item[ i)] $\left\|  b^{\ast}b\right\|  =\left\|  b\right\|  ^{2}$, and

\item[ j)] $b^{\ast}b\geq0\;\forall b\in$\emph{$\mathcal{B}$.}
\end{enumerate}
\end{definition}

Note that $\mathcal{B}_{0}$ is a $C^{\ast}$-algebra. Denote by $C_{c}%
(\mathcal{B})$ the collection of compactly supported continuous sections. Of
course, in our setting, because $\mathbb{Z}$ has the discrete topology,
continuous, compactly supported sections are really elements of the algebraic
direct sum $\sum\mathcal{B}_{k}$. Given $\xi,\eta\in C_{c}(\mathcal{B})$,
define the multiplication and involution by means of the formulae
\[
(\xi\ast\eta)(k)=\sum_{l}\xi(k-l)\eta(l),
\]%
\[
\xi^{\ast}(k)=\xi(-k)^{\ast}.
\]
Then $C_{c}(\mathcal{B})$ becomes a $\ast$-algebra. Let $P:C_{c}%
(\mathcal{B})\rightarrow\mathcal{B}_{0}$ be the restriction map, $\xi
\mapsto\xi(0)$. With the inner product $\langle\xi,\eta\rangle=P(\xi^{\ast
}\eta),\;C_{c}(\mathcal{B})$ becomes a pre-Hilbert $\mathcal{B}_{0}$-module.
For $\xi\in C_{c}(\mathcal{B})$, put $\left\|  \xi\right\|  _{2}:=\left\|
<\xi,\xi>\right\|  ^{1/2}$, and denote the completion of $C_{c}(\mathcal{B})$
with this norm by $L^{2}(\mathcal{B})$. Notice that
\[
L^{2}(\mathcal{B})=\bigoplus_{k\in\mathbb{Z}}\mathcal{B}_{k}%
\]
as Hilbert $\mathcal{B}_{0}$-modules. We have an embedding $C_{c}%
(\mathcal{B})\rightarrow L(L^{2}(\mathcal{B}))$ given by left multiplication.
Denote by $C^{\ast}(\mathcal{B})$ the completion of $C_{c}(\mathcal{B})$ with
respect to the operator norm. The map $P:C_{c}(\mathcal{B})\rightarrow
\mathcal{B}_{0},\;\xi\mapsto\xi(0)$ extends to a conditional expectation
$P:C^{\ast}(\mathcal{B})\rightarrow\mathcal{B}_{0}$.\medskip

\emph{Proof of the Theorem \ref{isomorphism}}. To the pair $(E_{\infty
},\mathcal{F}_{E})$, we can associate the Fell bundle $\mathcal{B}$, where
$\mathcal{B}_{n}:=E_{\infty}^{\otimes n},n\in\mathbb{Z}$. The multiplication
is given by the tensor product, where we identify $E_{\infty}^{\ast}\otimes
E_{\infty}$ with $\mathcal{F}_{E}$ and $E_{\infty}\otimes E_{\infty}^{\ast}$
with the ideal $\mathcal{F}_{E}^{1}$ of $\mathcal{F}_{E}$ generated by the
(images of) $K(E^{\otimes n})$, $n\geq1,$ in $\lim L(E^{\otimes n})$. (See
\cite{Pi}.) The involution is obvious. Then
\[
L^{2}(\mathcal{B})=\mathcal{E}_{\infty}=\bigoplus_{n\in\mathbb{Z}}E_{\infty
}^{\otimes n}.
\]
Since $\mathcal{E}_{\infty}$ is generated by $\mathcal{F}_{E}$ and $E_{\infty
}$, it follows that the $C^{\ast}$-algebra generated by the operators $M_{\xi
}$ is isomorphic to $C^{\ast}(\mathcal{B})$. Hence, $\tilde\mathcal{O}_{E}\simeq
C^{\ast}(\mathcal{B})$.

For the groupoid $\Gamma=\Gamma(X,T)$ and $l\in\mathbb{Z}$, take
\[
\Gamma_{l}:=\{(x,k,y)\in\Gamma\mid k=l\}=\{(x,y)\in X\times X\mid
x_{n}=y_{n+l}\;\mbox{for large}\;n\},
\]
and $\mathcal{D}_{l}=\overline{C_{c}(\Gamma_{-l})}$ (closure in $C^{\ast
}(\Gamma)$). Then it is easy to see that $C^{\ast}(\Gamma)$ is isomorphic to
$C^{\ast}(\mathcal{D})$. But $\mathcal{D}_{0}=C^{\ast}(R(X,T))\simeq\mathcal{F}_{E}=\mathcal{B}_{0}$, $\mathcal{D}_{1}=\overline
{C_{c}(\Gamma_{-1})}\simeq E_{\infty}=\mathcal{B}_{1}$, etc. Therefore the
Fell bundles $\mathcal{B}$ and $\mathcal{D}$ are isomorphic.

That concludes the proof of the theorem.
\endproof

\bigskip

Renault{\cite{Re2} has nice criteria for this $C^{\ast}$-algebra to be simple
and purely infinite. Applying them here, we obtain the following two propositions.

\begin{proposition}
Let $\sigma:X\rightarrow X$ be an essentially free branched covering. Assume
that for every nonempty open set $D\subset X$ and every $x\in X$, there exist
$m,n\in\mathbb{N}$ such that $x\in dom(T^{n})$ and $T^{n}x\in T^{m}D$. Then
$C^{\ast}(X,\sigma)$ is simple.
\end{proposition}

\begin{proposition}
Assume that for every nonempty open set $D\subset X$ there exists an open set
$D^{\prime}\subset D$ and $m,n\in\mathbb{N}$ such that  $T^{n}(\overline{D^{\prime}})$ is
strictly contained in $T^{m}(D^{\prime})$. Then $C^{\ast}(X,\sigma)$ is purely infinite.
\end{proposition}

The $C^{\ast}$-algebra of the equivalence relation $R(X,T)$ plays the role of
the AF-algebra. It is an inductive limit of unital non-continuous trace
algebras, similar to the dimension drop algebras considered by Elliott et al.

\section{K-theory computations and examples}

We use the six term exact sequence obtained by Pimsner to compute the K-theory
of the $C^{\ast}$-algebras associated to a branched covering.

\begin{theorem}
\label{P}Let $\sigma:X\rightarrow X$ be a branched covering of a locally
compact space $X$ with singular sets $S$ and $S^{\prime}$. Let $U=X\setminus
S$. We have an exact sequence
\[\begin{array}{ccccc}
K_{0}(C_{0}(U\cap\sigma(U)))&\longrightarrow &K_{0}(C_{0}(X))&\longrightarrow &
K_{0}(C^{\ast}(X,\sigma))\\
\uparrow &{}&{}&{}& \downarrow\\
K_{1}(C^{\ast}(X,\sigma))&\longleftarrow &K_{1}(C_{0}(X))&\longleftarrow
&K_{1}(C_{0}(U\cap\sigma(U)))\end{array}.
\]
\end{theorem}

\begin{proof}
Theorem 4.9 of \cite{Pi} for the augmented $C^{\ast}$-algebra $\tilde\mathcal{O}%
_{E}$ yields the following six term exact sequence
\[\begin{array}{ccccc}
K_{0}(I)& \overset{\otimes(\iota_{I}-[E])}{\longrightarrow}& K_{0}(C_{0}%
(X))& \overset{i_{\ast}}{\longrightarrow}& K_{0}(C^{\ast}(X,\sigma))\\
\delta\uparrow&{}&{}&{}& \downarrow\delta\\
K_{1}(C^{\ast}(X,\sigma))&\overset{i_{\ast}}{\longleftarrow}&K_{1}%
(C_{0}(X))&\overset{\otimes(\iota_{I}-[E])}{\longleftarrow}&K_{1}(I)\end{array}
\]
where $I$ is the ideal $I=\varphi^{-1}(K(E))\cap C_{0}(\sigma(U))$ in
$A=C_{0}(X)$, and where, recall, $C^{\ast}(X,\sigma)\simeq\tilde\mathcal{O}_{E}$.
The maps $\delta$ are the boundary maps corresponding to the Toeplitz
extension for $E$ described in \cite{Pi}. The map $i_{\ast}$ is the map on
$K$-theory induced by the inclusion of $C_{0}(X)$ in $C^{\ast}(X,\sigma)$, and
the two instances of $\otimes(\iota_{I}-[E])$ are the maps induced by Kasparov
multiplication by $\iota_{I}-[E]\in KK(I,C_{0}(X))$. (Here, $K_{\ast}%
(C_{0}(X))$ is identified with $KK_{\ast}(\mathbb{C},C_{0}(X))$, $K_{\ast
}(\mathcal{O}_{E})$ is identified with $KK_{\ast}(\mathbb{C},\mathcal{O}_{E}%
)$, and $K_{\ast}(I)$ is identified with $KK_{\ast}(\mathbb{C},I)$.) But
$\varphi^{-1}(K(E))=C_{0}(U)$, therefore $I=C_{0}(U\cap\sigma(U))$.
\end{proof}

\begin{example}
Consider the folding of the interval $\sigma:[0,1]\rightarrow [0,1]$.
Then \[C^*(R_0)=C([0,1]), C^*(R_1)=\{f:[0,1]\rightarrow {\bf M}_2\;\mid\; f(\frac{1}{2})\in {\mathbb C}\otimes I_2\},\]\[ C^*(R_2)=\{f:[0,1]\rightarrow {\bf M}_4\;\mid\; f(0), f(\frac{1}{4}), f(\frac{3}{4}), f(1) \in {\bf M}_2\otimes I_2, \; f(\frac{1}{2})\in {\mathbb C}\otimes I_4\}\] etc.

 Since $([0,1],\sigma)$ is essentially free, we may consider $C^{\ast
}([0,1],\sigma)$ defined above. Note that this $C^{\ast}$-algebra is not
simple, since the orbit of $0$ is $\{0,1\}$ which is not dense. To compute its
K-theory, note that $U=[0,1/2)\cup(1/2,1],\sigma(U)=[0,1)$, therefore
$I=C_{0}([0,1/2)\cup(1/2,1))$. The exact sequence becomes
\[\begin{array}{ccccc}
0 & \rightarrow & \mathbb{Z} & \rightarrow & K_{0}\\
\uparrow &{}&{}&{}& \downarrow\\
K_{1}& \leftarrow & 0 & \leftarrow & \mathbb{Z}\end{array},
\]
hence $C^{\ast}([0,1],\sigma)$ has $K_{0}=\mathbb{Z}\oplus\mathbb{Z}$ and
$K_{1}=0$.
\end{example}

\begin{example}
Let \[q(z)=\frac{(z^2+1)^2}{4z(z^{2}-1)},\;\; z\in{\bf S}^{2}={\mathbb C}\cup \{\infty\}.\] Then $q$ is a rational map of
degree $4$, which is a local homeomorphism, except at the points $z$ such that
the equation $q(z)=w$ has double roots. We find
\[
U={\bf S}^{2}\setminus\{\pm i, \pm(\sqrt{2}\pm 1)\},\;\;q(U)={\bf S}
^{2}\setminus\{ -1, 0, 1 \}.
\]
It is known (Latt\`es 1918) that the Julia set of $q$ is the whole Riemann sphere, and that every backward orbit $\displaystyle\bigcup_{n\geq 0}q^{-n}(z)$ is dense. From Theorem 4.2.5 in \cite{Be}, it follows that for any nonempty open set $W$ of ${\bf S}^2$, there is $N\geq 0$ such that $q^N(W)={\bf S}^2$. Notice that the forward orbits of all the points in the singular set $\{\pm i, \pm(\sqrt{2}\pm 1)\}$ are finite. Denote by $F$ the union of these orbits. Given any open set $D$,  we can find an open subset $D'\subset D$ such that the closure $\overline{D'}$ does not intersect the finite set $F$. We can find a positive integer $M$ such that  $q^M\mid_{D'}$ is a local homeomorphism and $q^M(D')={\bf S}^2 $, hence $\overline{D'}$ is strictly contained in $q^M(D')$.
In particular, $({\bf S}^2,q)$ satisfies the hypotheses of Propositions 3.6 and 3.7. Hence $C^*({\bf S}^2,q)$ is simple and purely infinite.
We have \[I \simeq C_{0}({\bf S}^{2}\setminus
\{\pm i, \pm(\sqrt{2}\pm 1), 0, \pm 1\}).\]
To compute the K-theory of $I$, we use the more general short exact sequence
\[
0\rightarrow C_{0}({\bf S}^{2}\setminus\{p_{1},p_{2},...,p_{n}
\})\rightarrow C({\bf S}^{2})\stackrel{j}{\rightarrow}{\mathbb C}^{n}\rightarrow 0.
\]
That gives 
\[\begin{array}{ccccc}
K_{0}&\rightarrow &{\mathbb Z}^{2}& \stackrel{j_*}{\rightarrow}& {\mathbb Z}^{n}\\\uparrow &{}&{}&{}& \downarrow\\
0 &\leftarrow &0 &\leftarrow &K_{1}\end{array}
\]
Using the fact that $j_*$ takes the Bott element into $0$, it follows that \[ker j_*\simeq {\mathbb Z}, \; coker j_*\simeq{\mathbb Z}^{n-1}.\] We obtain \[K_{0}(C_0({\bf S}^2\setminus\{p_1,p_2,...,p_n\}))={\mathbb Z},\; K_{1}(C_0({\bf S}^2\setminus\{p_1,p_2,...,p_n\}))={\mathbb Z}^{n-1}.\]
 The exact
sequence for our C*-algebra becomes 
\[\begin{array}{ccccc}
{\mathbb Z}&\rightarrow &{\mathbb Z}^{2}&\rightarrow &K_{0}(C^{\ast}({\bf S}
^{2},q))\\
\uparrow & {} &{} & {}& \downarrow \\
K_{1}(C^{\ast}({\bf S}^{2},q))& \leftarrow & 0 & \leftarrow & {\mathbb Z}^{8}\end{array}
\]

Note that for a general rational map $q$ of degree $\geq 2$, the Julia set may not be  the entire sphere. Its complement, the Fatou set, will provide an ideal in $C^*({\bf S}^2,q)$ such that the quotient will be a simple C*-algebra.
 For more details about the dynamics of a rational map, we refer to \cite{Be}.

\end{example}

\begin{example}
Consider $\mathbf{S}^{1}=\mathbb{R}\cup\{\infty\}$, and let $\sigma: \mathbb{R}\rightarrow
\mathbf{S}^{1},\sigma(t)=exp(2\pi it)$. Then $I=C_{0}(\mathbb{R})$ with
$K_{0}=0$ and $K_{1}=\mathbb{Z}$, therefore we have the exact sequence
\[\begin{array}{ccccc}
0&  \rightarrow & \mathbb{Z} & \rightarrow & K_{0}(C^{\ast}(\mathbf{S}^{1},\sigma))\\
\uparrow & {}&{}&{}& \downarrow\\
K_{1}(C^{\ast}(\mathbf{S}^{1},\sigma))& \leftarrow & \mathbb{Z} & \leftarrow
& \mathbb{Z}\end{array}
\]
Since the map $\mathbb{Z}\rightarrow \mathbb{Z}$ in the exact sequence is $id-id$, hence the zero map, it follows that
\[K_0(C^*({\bf S}^1, \sigma))\simeq \mathbb{Z}^2, \;\; K_1(C^*({\bf S}^1, \sigma))\simeq \mathbb{Z}.\]
Note that $C^*({\bf S}^1, \sigma)$ is simple, since every orbit is dense.
\end{example}

\end{document}